\documentclass[a4paper,12pt]{article}
\usepackage{amssymb,amsthm}
\usepackage{amsmath}

\def\<{\langle}
\def\>{\rangle}

\newcommand{\beq}{\begin{eqnarray}}
\newcommand{\eeq}{\end{eqnarray}}
\newcommand{\R}{\mathbb{R}}
\newcommand{\N}{\mathbb{N}}

\newtheorem{theorem}{\bf Theorem}
\newtheorem{lemma}[theorem]{\bf Lemma}

\newtheorem{corollary}[theorem]{\bf Corollary}
\newtheorem{remark}[theorem]{\bf Remark}
\newtheorem{example}[theorem]{\bf Example}

\def\|{\Vert}

\def\le {\leqslant}
\def\ge {\geqslant}

\def\part{\partial}

\begin{document}

\bigskip \bigskip \bigskip\bigskip\bigskip\bigskip

\title{Heat content, torsional rigidity and generalized Hardy inequalities
for complete Riemannian manifolds}

\author{Michiel van den Berg\\ \small
School of Mathematics,
  University of Bristol, \\\small University Walk,
\small Bristol BS8 1TW,
U.K.\\
\small  M.vandenBerg@bris.ac.uk\\\\
Peter B. Gilkey\\ \small Department of Mathematics,
  University of Oregon, \\
  \small Eugene, Oregon 97403,
  U.S.A.\\
\small gilkey@darkwing.uoregon.edu}

\date{}
\maketitle\vspace{10mm}

\normalsize

\begin{abstract}
We obtain upper bounds on the heat content and on the torsional
rigidity of a complete Riemannian manifold $M$, assuming a
generalized Hardy inequality for the Dirichlet Laplacian on
$M$.\\\\ 2000 Mathematics Subject Classification 35K05, 58G25,
60J65.
\end{abstract}

\begin{center}{\bf{Acknowledgements}}\end{center}
\small The first author was supported in part by the ESF. The
second author was supported in part by the NSF (U.S.A.). Both
authors wish to thank the Mittag-Leffler Institute for its
hospitality in the fall of 2002. The first author wishes thank
E.B. Davies for some useful discussions.

\normalsize
\section{Introduction}\label{sec1}
Let $M$ be a complete $C^\infty$ Riemannian manifold with boundary
$\partial M$, and let $-\Delta$ be the Dirichlet Laplacian on $M$.
In this paper we investigate the heat content of $M$ if $M$ has
uniform initial temperature one, while $\partial M$ is kept at
temperature zero for all time $t$. Let $u: M\times [0, \infty)\to
\R$ be the unique weak solution of

\begin{equation}\label{eq1}
\Delta u=\frac{\partial u}{\partial t},\qquad x\in M,\  t>0,
\end{equation}
\begin{equation}\label{eq2}
u(x,0)=1,\quad x\in M.
\end{equation}
The heat content $Q_M(t)$ is defined by
\begin{equation}\label{eq3}
Q_M(t)=\int_Mu(x;t)dx.
\end{equation}

It is well known that if $M$ is compact and $\partial M$ is
$C^\infty$ then there exists an asymptotic series for the heat
content of the form
\begin{equation}\label{eq4}
  Q_M(t)=\sum^J_{j=0} b_j t^{j/2}+O(t^{(J+1)/2}),\  t\to 0,
\end{equation}
where $J\in\N$, and where the coefficients are locally computable
invariants of $M$. In particular,
\begin{equation}\label{eq5}
  b_0=\int_M 1\  dx=\mbox{vol}(M),
\end{equation}
\begin{equation}\label{eq6}
b_1=-\frac{2}{\pi^{1/2}}\int_{\partial M}1\
dy=-\frac{2}{\pi^{1/2}}\mbox{area}(\partial M).
\end{equation}
For details we refer to \cite{4,5}.

In this paper we are concerned with the non-classical situation,
where $M$ is complete but non-compact. We shall be concerned with
the setting where either $M$ itself has infinite volume or where
$M$ has finite volume and $\partial M$ has infinite area. We
recall the following \cite{3}.

Let $M$ be the closure of an open set $M_0$ in Euclidean space
with boundary $\partial M=M\backslash M_0$. For $x\in M_0$ we
define the distance in direction $u, |u|=1$, by
\begin{equation}\label{eq7}
d_u(x)=\min\{|t|:x+tu\in \partial M\},
\end{equation}
and the mean distance function $\rho:M_0\to [0,\infty)$ by
\begin{equation}\label{eq8}
\frac{1}{\rho^2(x)}=\frac{1}{\mbox{area}(S_{m-1})}\int_{S_{m-1}}\frac{du}{d^2_u(x)},
\end{equation}
where $S_{m-1}$ is the unit sphere in $\R^m$. Theorem 3.4 in
\cite{3} asserts that if
\begin{equation}\label{eq9}
\rho^2(x)=o(\log(1+|x|^2))^{-1}, |x|\to\infty, x\in M_0
\end{equation}
then $M$ has finite heat content for all positive $t$. However, no
upper bounds for the heat content in terms of $t$ were obtained in
this general situation.

In this paper we obtain bounds for the heat content for a wide
class of $C^\infty$ complete Riemannian manifolds assuming a
generalized Hardy inequality.

Let $\delta :M\to[0,\infty)$ denote the distance to the boundary
function on $M$:
\begin{equation}\label{eq10}
  \delta(x)=\min\{d(x,y):y\in\partial M\},
\end{equation}
where $d(x,y)$ is the Riemannian distance between $x$ and $y$. We
say that $-\Delta$ satisfies a generalized Hardy inequality if
there exist constants $c>0$ and $\gamma\in (0,2]$ such that
\begin{equation}\label{eq11}
  -\Delta\ge \frac{c}{\delta^\gamma},
\end{equation}
in the sense of quadratic forms.
\begin{theorem}\label{The1}
Let $M$ be a complete $C^\infty$ Riemannian manifold, and suppose
that \eqref{eq11} holds for some $c>0$ and $\gamma\in (0,2]$.
Suppose there exists $\beta\in (0,2\gamma]$ such that
\begin{equation}\label{eq12}
\int_M\delta^\beta(x)dx<\infty.
\end{equation}
Then for all $t>0$
\begin{equation}\label{eq13}
  Q_M(t)\le \left(\frac{(\beta+\gamma)^2}{2e\beta\gamma
  c}\right)^{\beta/\gamma}\left(\int_M \delta^\beta(x)dx
  \right)t^{-\beta/\gamma}.
\end{equation}
\end{theorem}

\begin{theorem}\label{The2}
Let $M$ be a complete $C^\infty$ Riemannian manifold with finite
volume, and suppose that \eqref {eq11} holds for some $c>0$ and
$\gamma\in(0,2]$. Then for all $t>0$
\begin{equation}\label{eq14}
Q_M(t)\le {\rm{vol}}(M)-4^{-1}\int_{\{x\in M:\delta (x)<(2ct)^{1/\gamma}\}}
1\  dx.
\end{equation}
\end{theorem}
The Hardy inequality \eqref{eq11} guarantees that the boundary
$\partial M$ is not too thin, and that sufficient cooling of $M_0$
near $\partial M$ will take place. Condition \eqref{eq12} in
Theorem \ref{The1} guarantees that $M$ does not have to much
measure away from $\partial M$. Both the validity and applications
of inequalities like \eqref{eq11} to spectral theory have been
investigated in depth [8-11,\,13].

\begin{remark}\label{r3}
If $M_0$ is simply connected in $\R^2$ then \eqref{eq11} holds
with $\gamma=2$ and $c=1/16$. If $M_0$ is convex in $\R^m$ then
\eqref{eq11} holds with $\gamma=2$ and $c=1/4$ {\rm{\cite{6,7}}}.
\end{remark}
For open sets $M_0$ in $\R^m$ it was shown (Theorem
1.5.3 in \cite{6}) that
\begin{equation}\label{eq15}
  -\Delta\ge\frac{m}{4\rho^2},
\end{equation}
in the sense of quadratic forms. The proofs of Theorems \ref{The1}
and \ref{The2} together with \eqref{eq15} give the following.

\begin{corollary}\label{C3}
Let $M$ be the closure of an open set $M_0$ in Euclidean space
with boundary $\partial M = M\backslash M_0$. Suppose there exists
$\beta\in (0,4]$ such that
\begin{equation}\label{eq16}
\int_M \rho^\beta(x) dx<\infty.
\end{equation}
Then for all $t>0$
\begin{equation}\label{17}
  Q_M(t)\le\left(\frac{(\beta+2)^2}{e\beta
m}\right)^{\beta/2}\left(\int_M\rho^\beta(x) dx\right)t^{-\beta/2}.
\end{equation}
\end{corollary}

\begin{corollary}\label{C4}
Let $M$ be the closure of an open set $M_0$ in Euclidean space
with boundary $\partial M= M\backslash M_0$, and with finite
volume. Then for all $t>0$
\begin{equation}\label{18}
Q_M(t)\le {\rm{vol}}\ (M)-4^{-1} \int_{ \{x\in M:\rho(x)<(mt/2)^{1/2}\}}1\
dx,
\end{equation}
\end{corollary}
Let $p_M(x,y;t), x\in M, y\in M, t>0$  be the Dirichlet heat
kernel for $M$. We say that $M$ has finite torsional rigidity
$P_M$ if
\begin{equation}\label{eq19}
P_M=\int^\infty_0\int_M\int_M p_M(x,y;t)dx\; dy\;dt<\infty.
\end{equation}
Let $M_0$ be an open subset of $\R^m$. It is well known that if
$M_0$ has finite  volume then $M$ has finite torsional rigidity.
However, the converse is not true. In \cite{1} we showed that if
$M_0\subset\R^m$ satisfies \eqref{eq11} for $\gamma=2$ and some
$c>0$, then $P_M$ is finite if and only if \eqref{eq12} holds with
$\beta=2$. Since the solution of (\ref{eq1}-\ref{eq2}) with
Dirichlet boundary conditions on $\partial M$ is given by
\begin{equation}\label{eq20}
  u(x;t)= \int_Mp_M (x,y;t)dy,
\end{equation}
it follows that
\begin{equation}\label{eq21}
Q_M(t)=\int_M\int_M p_M(x,y;t)dx dy,
\end{equation}
and
\begin{equation}\label{eq22}
P_M=\int^\infty_0Q_M(t)dt.
\end{equation}
Theorem \ref{The1} gives the following.

\begin{corollary}\label{C5}
Suppose $M$ is a complete $C^\infty$ Riemannian manifold, and
suppose that \eqref{eq11} holds for some $c>0$ and
$\gamma\in(0,2]$. Suppose there exists $\varepsilon>0$ such that
\eqref{eq12} holds for all $\beta\in(\gamma-\varepsilon
,\gamma+\varepsilon)$.  Then $M$ has finite torsional rigidity.
\end{corollary}
The proof of Corollary \ref{C5} is elementary. We bound $Q_M(t)$
for small $t$ by \eqref{eq13} with $\beta=\gamma-\varepsilon/2$,
and for large $t$ by \eqref{eq13} with
$\beta=\gamma+\varepsilon/2$. Similarly one can show that for open
sets in $\R^m, P_M$ is finite if \eqref{eq16} holds for all
$\beta$ in some neighbourhood $\beta=2$. The following result is
an improvement.

\begin{theorem}\label{The6}
let $M$ be the closure of an open set $M_0$ in $\R^m$ with
boundary $\partial M=M\backslash M_0$. Suppose that
$\int_M\rho^2(x) dx<\infty$. Then
\begin{equation}\label{eq23}
  P_M\le \frac{4}{m}\int_M\rho^2(x) dx.
\end{equation}
\end{theorem}
In Lemma 2.6 of \cite{3} it was shown that if $M_0\subset\R^m$ is
such that $Q_M(t)$ is finite for all $t>0$ then trace
$(e^{t\Delta})$ is finite for all $t>0$. In the more general
setting of complete $C^\infty$ Riemannian manifolds we have the
following.

\begin{corollary}\label{C7}
Let $M$ be a complete $C^\infty$ Riemannian manifold, and suppose
that \eqref{eq11} holds for some $c>0$ and $\gamma\in (0,2].$
Suppose \eqref{eq12} holds for some $\beta\in(0,2\gamma]$, and
suppose there exists a function $g:(0,\infty)\to(0,\infty)$ such
that
\begin{equation}\label{eq24}
p_M(x,x;t)\le g(t), \qquad x\in M.
\end{equation}
Then for all $t>0$
\begin{equation}\label{eq25}
{\rm{trace}}(e^{t\Delta})=\int_M p_M(x,x;t)dx\le g(t/2)Q_M(t/2)<\infty.
\end{equation}
\end{corollary}
Sufficient conditions on the geometry of $M$ which guarantee the
uniform bound \eqref{eq24} were obtained by several authors
(Section 6 in \cite{10} and the references therein).

We conclude this introduction with an example to show that
Theorems \ref{The1} and \ref{The2} are close to being sharp.
\begin{example} Let $M(\alpha)\subset \R^2$ be given by
\begin{equation}\label{eq26}
  M(\alpha)=\{(\xi_1,\xi_2):\xi_1\ge 0, |\xi_2|\le(\xi_1+1)^{-\alpha}\},
\end{equation}
where $ \alpha>0$ is a constant. \end{example}
Since $M(\alpha)$ is
simply connected, we have by Remark \ref{r3} that \eqref{eq11}
holds with $\gamma=2$ and $c=\frac{1}{16}$. $M(\alpha)$ has
infinite volume if and only if $\alpha\le 1$. Estimate
\eqref{eq12} holds in this case if and only if
$\beta>\frac{1-\alpha}{\alpha}$. We can choose $\beta\in(0,4]$ if
and only if $\alpha>\frac{1}{5}$. We conclude by Theorem
\ref{The1} that for $\frac{1}{5}<\alpha\le 1$ and any
$\varepsilon>0$
\begin{equation}\label{eq27}
  Q_{M(\alpha)}(t)\le K_1
  t^{(\alpha-1)/(2\alpha)-\varepsilon},\quad t>0,
\end{equation}
where $K_1$ is a finite positive constant depending on $\alpha$
and $\varepsilon $ respectively.

Theorem \ref{The2} gives that for $\alpha>1$
\begin{equation}\label{eq28}
  \liminf_{t\to
  0}(\mbox{vol}(M(\alpha)-Q_{M(\alpha)}(t))t^{(1-\alpha)/(2\alpha)}\ge
  K_2,
\end{equation}
where $K_2$ is a strictly positive constant depending on $\alpha$.
The precise asymptotic behaviour of $Q_{M(\alpha)}(t)$ as $t\to 0$
has been computed in \cite{2}. The results in \cite{2} show for
example that \eqref{eq27} holds for all $ 0<\alpha<1$ with
$\varepsilon=0$.

\section{Proof of Theorem 1}\label{Sec2}

Let $\Delta$ be the Dirichlet Laplacian acting in $L^2(M)$, and
let $u:M\times [0,\infty)\to\R$ be the unique weak solution
of\eqref{eq1} with initial condition
\begin{equation}\label{eq29}
u(x;0)=f(x),\qquad x\in M,
\end{equation}
where $f:M\to[0,\infty)$ is bounded and measurable. Then
\begin{equation}\label{eq30}
  u=e^{t\Delta}f,
\end{equation}
with integral representation
\begin{equation}\label{eq31}
  u(x;t)=\int_Mp_M(x,y;t)f(y)dy.
\end{equation}
We let $\varepsilon>0$, and choose
$\{f_\varepsilon:\varepsilon>0\}$ to be a family of $C^\infty$
functions on $M$ such that $0\leq f_\varepsilon\leq 1,\
f_\varepsilon$ is monotone increasing as $\varepsilon\to 0$, and
\begin{equation}\label{eq32}
f_\varepsilon(x)=\left\{\begin{array}{ll}
0, & \delta(x)<\varepsilon,\\
1,&2\varepsilon\leq \delta(x).

\end{array}\right.
\end{equation}
It follows by the maximum principle (Section 2.4 in \cite{10})
that
\begin{equation}\label{eq33}
  \int_Mp_M(x,y;t)dy\leq 1.
\end{equation}
By Fubini's theorem and \eqref{eq33} we conclude that the unique
solution $u_\varepsilon$ satisfies
\begin{align}\label{eq34}
\begin{split}
\int_M u_\varepsilon(x;t)dx&=\int_M\int_M
p(x,y;t)f_\varepsilon(y)\; dx\; dy\\ &\leq \int_M
f_\varepsilon(y)dy\leq\int_{\{x\in M:
\delta(x)\geq\varepsilon\}}1\  dx\\ &\leq \int_{{\{x\in M:
\delta^\beta
(x)\ge\varepsilon^\beta\}}}\varepsilon^{-\beta}\delta^\beta
(x)\leq\varepsilon^{-\beta} \int_M\delta^\beta(x) dx<\infty.
\end{split}
\end{align}
Let $p\ge 3/2$. By the maximum principle we have that $0\leq
u_\varepsilon\leq 1$. Hence $u^{2p-2}_\varepsilon\leq
u_\varepsilon$, and
\begin{equation}\label{eq35}
\int_M u_\varepsilon^{2p-2}(x;t) dx\leq \int_M
u_\varepsilon(x;t)dx.
\end{equation}

Hence, by Fubini's theorem, Cauchy-Schwarz's inequality and
\eqref{eq35}
\begin{align}\label{eq36}
\begin{split}
\left|-\frac{d}{dt}\int_M
u^p_\varepsilon(x;t)dx\;\right|&=p\left|\int_M
u^{p-1}_\varepsilon(x;t)\frac{\partial u_\varepsilon}{\partial
t}(x;t)dx\right|\\ &\leq p\left\{\int_M
u_\varepsilon^{2p-2}(x;t)dx\right\}^{1/2}\left\{\int_M\left(\frac{\partial
u_\varepsilon}{\partial t}(x;t)\right)^2dx\right\}^{1/2}\\ & \leq
p\left\{\int_Mu_\varepsilon(x;t)dx\right\}^{1/2}\left\{\int_M\left(\frac{\partial
u_\varepsilon}{\partial t}(x;t)\right)^2 dx\right\}^{1/2}.
\end{split}
\end{align}
But
\begin{align}\label{eq37}
\begin{split}
\int_M\left(\frac{\partial u_\varepsilon}{\partial
t}(x;t)\right)^2 dx &=\langle\Delta e^{t\Delta}f_\varepsilon,\;
\Delta e^{t\Delta}f_\varepsilon\rangle\\ &=\langle(\Delta
e^{t\Delta/2})e^{t\Delta/2}f_\varepsilon, (\Delta
e^{t\Delta/2})e^{t\Delta/2}f_\varepsilon\rangle,
\end{split}
\end{align}
and since $\Delta e^{t\Delta/2}$ is a bounded operator, bounded by
$2/(et)$, we have that
\begin{align}\label{eq38}
\begin{split}
\int_M \left( \frac{\partial u_\varepsilon}{\partial t}(x;t)
\right)^2 dx &\leq \frac{4}{e^2t^2}\langle
e^{t\Delta/2}f_\varepsilon, e^{t\Delta/2} f_\varepsilon\rangle\\
&\leq \frac{4}{e^2t^2}\langle e^{t\Delta}f_\varepsilon,
f_\varepsilon\rangle\\ &\leq \frac{4}{e^2t^2}\langle
e^{t\Delta}f_\varepsilon, 1\rangle\\ &=\frac{4}{e^2t^2}\int_M
u_\varepsilon(x;t)dx.
\end{split}
\end{align}
By (\ref{eq36}-\ref{eq38}) we conclude that we have the estimate
\begin{equation}\label{eq39}
-\frac{d}{dt}\int_Mu_\varepsilon^p(x;t)dx\leq\frac{2p}{et}\int_M
u_\varepsilon(x;t)dx.
\end{equation}

We use Fubini's theorem, integration by parts, and the generalized
Hardy inequality \eqref{eq11} to obtain
\begin{align}\label{eq40}
\begin{split}
-\frac{d}{dt}\int_M u_\varepsilon^p(x;t)dx &= -p\int_M
u^{p-1}_\varepsilon(x;t)\Delta u_\varepsilon(x;t) dx\\ &= p\;(p-1)
\int_Mu_\varepsilon^{p-2}(x;t)(\nabla u_\varepsilon(x;t))^2 dx\\
&=\frac{4(p-1)}{p}\int_M|\nabla u_\varepsilon^{p/2}(x;t)|^2 dx\\
&\geq \frac{4(p-1)c}{p}\int_M
u^p_\varepsilon(x;t)\delta^{-\gamma}(x) dx.
\end{split}
\end{align}
Combining the estimates of \eqref{eq39} and \eqref{eq40} yields
\begin{equation}\label{eq41}
  \int_M u_\varepsilon(x;t)\geq \frac{2(p-1) ect}{p^2}\int_M
  u_\varepsilon^p(x;t)\delta^{-\gamma}(x)dx.
\end{equation}
By \eqref{eq34} the left hand side of \eqref{eq41} is finite and
this implies that the right hand side of \eqref{eq41} is finite as
well. By H\"{o}lder's inequality
\begin{equation}\label{eq42}
\int_M u_\varepsilon(x;t)dx \leq \left\{\int_M
u_\varepsilon^p(x;t)\delta^{-\gamma}(x)
dx\right\}^{1/p}\left\{\int_M \delta(x)^{\gamma/(p-1)}dx\right\}^{(p-1)/p}.
\end{equation}
Since $\beta \in (0,2\gamma]$ by the hypothesis in Theorem
\ref{The1} the choice
\begin{equation}\label{eq43}
p=1+\frac{\gamma}{\beta}
\end{equation}
guarantees that $p\ge3/2$. Then $\delta^{\gamma/(p-1)}$ is
integrable by \eqref{eq12}. By \eqref{eq41} and \eqref{eq42}
\begin{equation}\label{eq44}
\int_M u^p_\varepsilon(x;t)\delta^{-\gamma}(x)
dx\le\left(\frac{p^2}{2(p-1)ect}\right)^{p/(p-1)}\int_M\delta (x)
^{\gamma/(p-1)} dx ).
\end{equation}
Substitution of \eqref{eq44} into \eqref{eq42} then results into
\begin{equation}\label{eq45}
\int_M u_\varepsilon(x;t) dx \leq \left(\frac{p^2}{2(p-1)
ect}\right)^{1/(p-1)}\int_M \delta(x)^{\gamma/(p-1)}dx,
\end{equation}
and a further substitution of \eqref{eq43} into \eqref{eq45} gives
\begin{equation}\label{eq46}
  \int_M u_\varepsilon(x;t) dx\leq
  \left(\frac{(\beta+\gamma)^2}{2e\beta\gamma c
t}\right)^{\beta/\gamma}\int_M\delta^\beta(x)
  dx.
\end{equation}

Since the right hand side of \eqref{eq46} is independent of
$\varepsilon$ we have by Fatou's lemma
\begin{align}\label{eq47}
\begin{split}
Q_M (t)&= \int_M(\lim_{\varepsilon\to 0} u_\varepsilon(x;t))dx\\ &
\leq\liminf_{\varepsilon\to 0}\int_M u_\varepsilon(x;t)\leq
\left(\frac{(\beta+\gamma)^2}{2e\beta\gamma c
t}\right)^{\beta/\gamma}\int_M \delta^\beta(x) dx.
\end{split}
\end{align}
\qed
\section{Proof of Theorem 2}\label{Sec3}

Let $f=1$ in \eqref{eq29}, and note that
\begin{equation}\label{eq48}
Q_M(2t)=\langle e^{2t\Delta} 1, 1\rangle=\langle
e^{t\Delta}1,e^{t\Delta}1\rangle=\int_M u^2(x;t) dx.
\end{equation}
Hence \eqref{eq40} for $p=2$ gives
\begin{equation}\label{eq49}
  -\frac{d}{dt}Q_M(2t)\ge 2c\int_M u^2(x;t) \delta^{-\gamma}(x)
  dx.
\end{equation}
Integrating this inequality with respect to $t$ over $[0,t]$
yields by Fubini's theorem
\begin{equation}\label{eq50}
  {\rm{vol}} (M) -Q_M (2t)\ge 2c\int_M dx \int_0^t u^2
  (x;\tau)\delta^{-\gamma}(x) d\tau.
\end{equation}
It follows, by the maximum principle, that $u(x;\tau)\ge u(x;2t)$
for all $0\le \tau\le t$. Hence
\begin{equation}\label{eq51}
  {\rm{vol}} (M)-Q_M (2t)\ge 2 ct \int_M u^2
  (x;2t)\delta^{-\gamma}(x) dx.
\end{equation}
Let $u=1-v$. We use the definition given in \eqref{eq3}
and we use equation \eqref{eq51} to see that for any $\varepsilon>0$
\begin{align}\label{eq52}
\begin{split}
\int_Mv (x;t) dx & \ge c t \int_M(1-v(x;t))^2 \delta^{-\gamma}(x)
dx\\ & \geq ct \int_{\{x\in M: \delta(x)<\varepsilon\}}
(1-\upsilon(x;t))^2 \delta^{-\gamma}(x) dx\\ & \geq ct
\varepsilon^{-\gamma}\int_{\{x\in
M:\delta(x)<\varepsilon\}}(1-2\upsilon(x;t)) dx).
\end{split}
\end{align}
It follows that
\begin{equation}\label{eq53}
(1+2ct\varepsilon^{-\gamma})\int_M v(x;t) dx\ge c t
\varepsilon^{-\gamma}\int_{\{x\in M: \delta(x)<\varepsilon\}}1\
dx.
\end{equation}
The choice
\begin{equation}\label{eq54}
\varepsilon=(2ct)^{1/\gamma}
\end{equation}
in \eqref{eq53} completes the proof of Theorem \ref{The2}. \qed\\

To prove Corollaries \ref{C3} and \ref{C4} respectively we note
that \eqref{eq15} holds for open sets $M_0\subset\R^m$. We follow
the proofs of Theorems \ref{The1} and \ref{The2} respectively by
replacing $\delta$ by $\rho$, $c$ by $m/4$ and $\gamma$ by $2$
throughout.

To prove Corollary \ref{C7} we note that by (5.2) in Section 5.1
of \cite {10} and \eqref{eq24}
\begin{equation}\label{eq55}
p_M (x,y;t) \le (p_M(x,x;t)p_M(y,y;t))^{1/2}\le g(t).
\end{equation}
Hence
\begin{align}\label{eq56}
\begin{split}
\int_M p_M (x,x;t) dx &=\int_M\int_M p^2_M (x,y;t/2)dx dy\\&\leq
g(t/2)\int_M\int_M p_M (x,y;t/2) dx dy \leq g(t/2)Q_M(t/2)<\infty.
\end{split}
\end{align}
\qed
\section{Proof of Theorem \ref{The6}}\label{Sec4}
The proof of Theorem \ref{The6} is based on a couple of lemmas
which are of independent interest, and which are related to
results on the expected life time of $h$-conditioned Brownian
motion \cite{12,13}.

\begin{lemma}\label{L8}
Let $M_0$ be an open set in $\R^m, m\ge 3$. Suppose that for all
$\varepsilon>0$
\begin{equation}\label{eq57}
{\rm{vol}} (\{x\in M_0: \delta(x)>\varepsilon\})<\infty.
\end{equation}
Let $\{f_\varepsilon:\varepsilon>0\}$ be as in the proof of
Theorem \ref{The1}. Then
\begin{equation}\label{eq58}
-\Delta w=f_\varepsilon,
\end{equation}
has a unique, weak, bounded and non-negative solution
$w_\varepsilon$ with
\begin{equation}\label{eq59}
||w_\varepsilon||_\infty\le \frac{m}{4\pi(m-2)}{\rm{vol}} \{x\in
M:\,
\delta(x)>\varepsilon\}^{2/m}.
\end{equation}
\end{lemma}
\begin{proof}
Put
\begin{equation}\label{eq60}
M_\varepsilon=\{x\in M: \delta(x)>\varepsilon\}.
\end{equation}
Since $\delta\leq \rho$, \eqref{eq68} implies that
$\int\limits_M\delta^2(x)dx< \infty$. Then the last three
inequalities in the right hand side of \eqref{eq34} for $\beta=2$
imply ${\rm{vol}}(M_\varepsilon)<\infty$ for $\varepsilon>0$.

Since $(-\Delta)^{-1}$ has integral kernel
\begin{equation}\label{eq61}
\int_0^\infty p_M (x,y;t) dt
\end{equation}
we have that
\begin{align}\label{eq62}
\begin{split}
w_\varepsilon(x)&= \int_M\int_0^\infty p_M(x,y;t)\; dt\;
f_\varepsilon(y) dy\\ & = \int_{M_\varepsilon}\int^\infty_0
p_M(x,y;t)\;dt\; f_\varepsilon(y) dy\\ &\leq
\int_{M_\varepsilon}\int^\infty_0 p_M(x,y;t)\;dt\; dy.
\end{split}
\end{align}
By positivity of the Dirichlet heat Kernel
\begin{equation}\label{eq63}
  \int_{M_\varepsilon} p_M (x,y;t)dy \leq \int_M p_M (x,y;t)dy
  \leq 1.
\end{equation}
Hence for $t_0>0$ we have by Fubini's theorem
\begin{equation}\label{eq64}
\int_{M_\varepsilon}\int _0^{t_0}p_M (x,y;t)\;dt\; dy \leq t_0.
\end{equation}
Moreover, by monotonicity of the Dirichlet heat Kernel
\begin{equation}\label{eq65}
  p_M(x,y;t)\le (4\pi t)^{-m/2}.
\end{equation}
Hence for $m\geq 3$
\begin{equation}\label{eq66}
\int_{M_\varepsilon}\int^\infty_{t_0} p_M (x,y;t)\; dt\; dy \leq
(4\pi)^{-m/2}\frac{2}{m-2}t_0^{1-m/2}{\rm{vol}}(M_\varepsilon).
\end{equation}
By \eqref{eq64} and \eqref{eq66} we conclude that
\begin{equation}\label{eq67}
  w_\varepsilon(x)\leq
  t_0+(4\pi)^{-m/2}\frac{2}{m-2}t_0^{1-m/2}{\rm{vol}}(M_\varepsilon).
\end{equation}
We minimize the right hand side of \eqref{eq67} by setting
\begin{equation}\label{eq68a}
t_0=(4\pi)^{-1}({\rm{vol}}(M_\varepsilon))^{2/m},
\end{equation}
and \eqref{eq59} now follows.
\end{proof}

\begin{lemma}\label{L9}
Let $M_0$ be an open set in $\R^2$. Suppose that
\begin{equation}\label{eq68}
\int_M \rho^2(x) dx<\infty.
\end{equation}
Let $\{f_\varepsilon: \varepsilon>0\}$ be as in the Proof of
Theorem \ref{The1}. Then \eqref{eq58} has a unique, weak, bounded
and non-negative solution $w_\varepsilon$ with
\begin{equation}\label{eq69}
||w_\varepsilon||_\infty\leq
\left(\frac{8}{\pi}\,{{\rm{vol}}(M_\varepsilon)}\int_M\rho^2(x) dx
\right)^{1/2}.
\end{equation}
\end{lemma}
\begin{proof}
Let
\begin{equation}\label{eq70}
\lambda=\rm{inf\, spec} (-\Delta).
\end{equation}
Then by \eqref{eq15} we have for any smooth function $u$ with
compact support in $M_0$, and with $||u||_2=1$,
\begin{align}\label{eq71}
\begin{split}
1=||u||^2_2&\leq \left(\int_M
\left(\frac{u(x)}{\rho(x)}\right)^2dx \int_M\rho^2(x) dx
\right)^{1/2}\\ &\leq \left(2\int_M|\nabla u(x)|^2 dx \int_M
\rho^2 (x) dx \right)^{1/2}.
\end{split}
\end{align}
Taking the infimum over all such $u$ we obtain by \eqref{eq68} and
\eqref{eq71}
\begin{equation}\label{eq72}
\lambda\geq \frac{1}{2}\left(\int_M \rho^2(x)dx\right)^{-1}>0.
\end{equation}
To prove Lemma \ref{L9} we have by domain monotonicity
\begin{equation}\label{eq73}
p_M(x,x;t) \leq e^{-t\lambda/2} p_M(x,x;t/2)\leq (2\pi t)^{-1}
e^{-t\lambda/2},
\end{equation}
and hence by \eqref{eq55}
\begin{equation}\label{eq74}
p_M(x,y;t)\leq(2\pi t)^{-1} e^{-t\lambda/2}.
\end{equation}
As in the proof of Lemma \ref{L8} we have the estimates
(\ref{eq62}-\ref{eq64}). Estimate \eqref{eq66} is replaced, using
\eqref{eq74}, by

\begin{align}\label{eq75}
\begin{split}
\int_{M_\varepsilon} \int_{t_0}^\infty p_M(x,y;t) dt dy &\leq
\int_{M_\varepsilon}1\, dy \int_{t_0}^\infty (2\pi
t)^{-1}e^{-t\lambda/2}dt\\ &= (\pi
t_0\lambda)^{-1}{\rm{vol}}(M_\varepsilon).
\end{split}
\end{align}
By (\ref{eq62}-\ref{eq64}) and \eqref{eq75}
\begin{equation}\label{eq76}
w_\varepsilon(x)\leq
t_0+(\pi t_0\lambda)^{-1}{\rm{vol}}(M_\varepsilon).
\end{equation}
We minimize the right hand side of \eqref{eq76} by setting
\[t_0=({\rm{vol}}(M_\varepsilon)/(\pi\lambda))^{1/2}.\] We then
substitute the lower bound for $\lambda$ in \eqref{eq72} to obtain
\eqref{eq69}.
\end{proof}

We can now give the proof of Theorem \ref{The6}. We choose the
family \\ $\{f_\varepsilon: \varepsilon>0\}$ as in the proof of
the Theorem \ref{The1}. Since we have assumed that
$\int_M\rho^2(x) dx<\infty$ the support of $f_\varepsilon$ has
finite volume. Since $0\leq f_\varepsilon\leq 1$ and since
$||w_\varepsilon||_\infty<\infty$ (by Lemmas \ref{L8} and
\ref{L9}) we have that $w_\varepsilon f_\varepsilon$ is integrable
on $M$. By \eqref{eq58}, an integration by parts and \eqref{eq15}
we conclude that

\begin{align}\label{eq77}
\begin{split}
\infty> \int_M w_\varepsilon(x)f_\varepsilon(x) dx\geq
-\int_{M}w_\varepsilon(x)\Delta w_\varepsilon(x)
dx&=\int_{M} |\nabla w_\varepsilon(x)|^2 dx\\
&\geq
\frac{m}{4}\int_{M} w^2_\varepsilon(x)\rho^{-2} (x) dx.
\end{split}
\end{align}
By assumption $\int_M\rho^2(x)dx<\infty$. Hence by H\"{o}lder's
inequality

\begin{align}\label{eq78}
\begin{split}
\int_M w_\varepsilon(x)f_\varepsilon(x)dx&\leq
\left(\int_{M}\frac{w_\varepsilon^2(x)f^2_\varepsilon(x)}{\rho^2(x)}dx\int_M
\rho^2(x) dx\right)^{1/2}\\ &\leq
\left(\int_{M}w^2_\varepsilon(x)\rho^{-2}(x) dx \int_M \rho^2(x)
dx\right)^{1/2}
\end{split}
\end{align}
By \eqref{eq77} and \eqref{eq78}
\begin{equation}\label{eq79}
\int_M w_\varepsilon^2(x)\rho^{-2}(x) dx \leq
\frac{16}{m^2}\int_M \rho^2(x) dx,
\end{equation}
and by \eqref{eq78} and \eqref{eq79}
\begin{equation}\label{eq80}
\int_M w_\varepsilon(x)f_\varepsilon(x) dx \leq \frac{4}{m}\int_M
\rho^2(x)dx.
\end{equation}
By the first equality in \eqref{eq62}, \eqref{eq80} and Fubini's
theorem
\begin{equation}\label{eq81}
\int\limits_0^\infty \int\limits_M\int\limits_M p_M (x,y;t) f_\varepsilon(x)
f_\varepsilon(y)\;dx\;dy\;dt \leq
\frac{4}{m}\int\limits_M \rho^2(x) dx.
\end{equation}
Since $f_\varepsilon$ is increasing as $\varepsilon\to 0$, we have
by the monotone convergence theorem
\begin{equation}\label{eq82}
  \int\limits_0^\infty\int\limits_M\int\limits_M p_M (x,y;t)\;dx\;dy\;dt\le
  \frac{4}{m}\int\limits_M \rho^2(x) dx.
\end{equation}
But this is the conclusion of Theorem \ref{The6} by definition
\eqref{eq22}. \qed


\begin{thebibliography}{99}
\bibitem{1}{\sc R. Ba\~{n}uelos, M. van den Berg, T. Carroll}, \emph
{Torsional rigidity and expected lifetime of
Brownian motion}, J. London Math. Soc. (2) {\bf 66} (2002)
499--512.
\bibitem{2}{\sc M. van den Berg}, \emph {Heat content asymptotics for planar
regions with
cusps}, J. London Math. Soc. (2) {\bf 57} (1998) 677--693.
\bibitem{3}{\sc M. van den Berg, E.B. Davies}, \emph {Heat flow out of
regions in
$\R^m$}, Math. Zeit. {\bf 202} (1989) 463--482.
\bibitem{4}{\sc M. van den Berg, P.B. Gilkey}, \emph {Heat content
asymptotics of a Riemannian manifold with
boundary}, J. Funct. Anal. {\bf 120} (1994) 48--71.
\bibitem{5}{\sc M. van den Berg, P.B. Gilkey}, \emph {The heat equation with
inhomogeneous Dirichlet boundary
conditions}, Comm. Analysis and Geometry {\bf 7} (1999) 279--294.
\bibitem{12}{\sc M.Cranston, T.R. McConnel}, \emph {The lifetime of
conditioned Brownian
motion}, Z. Wahrsch. Verw. Gebiete {\bf 65} (1983) 1--11.
\bibitem{13}{\sc M.Cranston}, \emph {Lifetime of conditioned
Brownian motion in Lipschitz domains}, Z. Wahrsch. Verw. Gebiete
{\bf 70} (1985) 335--340.
\bibitem{6}{\sc E.B. Davies}, \emph {Heat Kernels and spectral
theory}. Cambridge Tracts in Mathematics, {\bf 92}. Cambridge
University Press, Cambridge, 1989.
\bibitem{7}{\sc E.B. Davies}, \emph {A review of Hardy
inequalities}, Operator Theory Adv. Appl. {\bf 110} (1999) 55--67.
\bibitem{8}{\sc E.B. Davies}, \emph {Sharp boundary estimates for
elliptic operators}, Math. Proc. Cambridge Philos. Soc. {\bf 129}
(2000) 165--178.
\bibitem{9}{\sc E.B. Davies, B. Simon}, \emph {Ultracontractivity
and the heat kernel for Schr\"{o}dinger operators and Dirichlet
Laplacians}, J. Funct. Anal. {\bf 59} (1984) 335--395.
\bibitem{10}{\sc A. Grigoryan}, \emph {Estimates of heat kernels
on Riemannian manifolds}, Spectral theory and geometry (Edinburgh
1998), 140--225, London Math. Soc. Lecture Note Series {\bf 273}.
Cambridge Univ. Press, Cambridge, 1999.
\bibitem{11}{\sc M. Lianantonakis}, \emph {On the eigenvalue
counting function for weighted Laplace-Beltrami operators}, J.
Geom. Anal. {\bf 10} (2000) 299--322. \nopagebreak
\bibitem{14} {\sc G. P\'{o}lya, G. Szeg\"{o}}, \emph {Isoperimetric
inequalities in mathematical
physics}. Princeton University Press, Princeton, 1951.

\end{thebibliography}
\end{document}